\documentclass[12pt]{article}
\usepackage{amsmath}
\usepackage{amssymb}
\usepackage{vatola}

\def\q{\quad}

\def\qtq#1{\q\t{#1}\q}
\def\mod#1{\ (\text{\rm mod}\ #1)}
\def\t{\text}
\def\f{\frac}
\def\e{\equiv}
\def\b{\binom}

\def\ap{\langle a\rangle_p}
\def\sls#1#2{(\f{#1}{#2})}
 
\def\Ls#1#2{\Big(\f{#1}{#2}\Big)}

\let \pro=\proclaim
\let \endpro=\endproclaim

\begin{document}
\par\q\par\q
 \centerline {\bf
Note on super congruences modulo $p^2$}
$$\q$$
\centerline{Zhi-Hong Sun} $$\q$$ \centerline{School of Mathematical
Sciences}
 \centerline{Huaiyin Normal University} \centerline{Huaian, Jiangsu 223001, P.R.
China} \centerline{Email: zhihongsun@yahoo.com}
\centerline{Homepage: http://www.hytc.edu.cn/xsjl/szh} \centerline{}

 \abstract{Let $p$ be an odd prime, and let $m$ be an integer
 with $p\nmid m$. In this paper show that
$$\sum_{k=0}^{p-1}\frac{\binom{2k}k\binom ak\binom{-1-a}k}{m^k}
\equiv 0\pmod p
 \quad\hbox{implies}\quad\sum_{k=0}^{p-1}\frac{\binom{2k}k\binom ak
 \binom{-1-a}k}{m^k}\equiv 0\pmod  {p^2}.$$
\par\q
 }
 \endabstract

\section*{1. Introduction}
\par\q  Let $[x]$ be the greatest integer not exceeding $x$.
 For a prime $p$ let $\Bbb
Z_p$ be the set of rational numbers whose denominator is not
divisible by $p$. Let $\{P_n(x)\}$ be the Legendre polynomials given
by
$$P_0(x)=1,\ P_1(x)=x,\ (n+1)P_{n+1}(x)=(2n+1)xP_n(x)-nP_{n-1}(x)\ (n\ge
1).$$ It is well known that (see [MOS, pp.\;228-232], [B1] and [B2])
$$P_n(x)=\f
1{2^n}\sum_{k=0}^{[n/2]}\b nk(-1)^k\b{2n-2k}nx^{n-2k}
=\sum_{k=0}^n\b nk\b{n+k}k\Ls{x-1}2^k.\tag 1.1$$ From (1.1) we see
that $P_n(-x)=(-1)^nP_n(x)$.
\par Let $p>3$ be a prime. In 2003  Rodriguez-Villegas [RV] conjectured that
$$\aligned&\sum_{k=0}^{p-1}\f{\b{2k}k^2\b{3k}k}{108^k}\e 0\mod{p^2}\qtq{for}
p\e 2\mod 3,
\\& \sum_{k=0}^{p-1}\f{\b{2k}k^2\b{4k}{2k}}{256^k} \e 0\mod{p^2}
\qtq{for} p\e 5,7\mod 8,
\\& \sum_{k=0}^{p-1} \f{\b{2k}k\b{3k}k\b{6k}{3k}}{1728^k}
\e 0\mod{p^2}\qtq{for} p\e 3\mod 4.\endaligned\tag 1.2$$ In 2005
Mortenson [M] proved these congruences modulo $p$, in 2012 Z.W. Sun
[Su3] confirmed (1.2). Motivated by Mortenson's work, in [Su1,Su2]
Z.W. Sun posed many conjectures concerning the following sums modulo
$p^2$:
$$\sum_{k=0}^{p-1}\f{\b{2k}k^3}{m^k},\
\sum_{k=0}^{p-1}\f{\b{2k}k\b{3k}k}{m^k},\
\sum_{k=0}^{p-1}\f{\b{2k}k\b{4k}{2k}}{m^k}\qtq{and} \sum_{k=0}^{p-1}
\f{\b{2k}k\b{3k}k\b{6k}{3k}}{m^k},$$ where $m$ is an integer with
$p\nmid m$. In [S2-S5] the author solved some of his conjectures by
establishing the following congruences:
$$\align &\sum_{k=0}^{p-1}\f{\b{2k}k^3}{m^k}
\e P_{\f{p-1}2}\Big(\sqrt{1-\f{64}m}\Big)^2\mod {p^2},\tag 1.3
\\&\sum_{k=0}^{p-1}\b{2k}k^2\b{3k}k(x(1-27x))^k\e\Big(
\sum_{k=0}^{p-1}\b{2k}k\b{3k}kx^k\Big)^2\mod{p^2},\tag 1.4
\\&\sum_{k=0}^{p-1}\b{2k}k^2\b{4k}{2k}(x(1-64x))^k\e\Big(
\sum_{k=0}^{p-1}\b{2k}k\b{4k}{2k}x^k\Big)^2\mod{p^2},\tag 1.5
\\&\sum_{k=0}^{p-1}\b{2k}k\b{3k}k\b{6k}{3k}(x(1-432x))^k\e\Big(
\sum_{k=0}^{p-1}\b{2k}k\b{3k}k\b{6k}{3k}x^k\Big)^2\mod{p^2}.\tag 1.6
\endalign$$
 It is easily seen that (see [S1, pp.1916-1917, 1920], [S3, p.1953] and [S5, p.182])
$$\aligned&\b{-\f 12}k^2=\f{\b{2k}k^2}{16^k},\ \b{-\f
13}k\b{-\f 23}k=\f{\b{2k}k\b{3k}k}{27^k}, \\&\ \b{-\f 14}k\b{-\f
34}k=\f{\b{2k}k\b{4k}{2k}}{64^k},\ \b{-\f 16}k\b{-\f
56}k=\f{\b{3k}k\b{6k}{3k}}{432^k}.\endaligned\tag 1.7$$ Let $p$ be
an odd prime, $a\in\Bbb Z_p$ and let $\ap\in\{0,1,\ldots,p-1\}$ be
given by $a\e \ap\mod p$. In [S6], the author proved that
$$\sum_{k=0}^{p-1}\f{\b{2k}k\b ak\b{-1-a}k}{4^k}\e 0\mod{p^2}
\qtq{for}\ap\e 1\mod 2.$$
 Let
$(a)_k=a(a+1)\cdots(a+k-1)=(-1)^kk!\b{-a}k$ and
$${}_rF_s\biggl(\begin{matrix} a_1,\ldots,a_r\\
b_1,\ldots,b_s\end{matrix}
\biggm|z\biggr)=1+\sum_{k=1}^{\infty}\f{(a_1)_k\cdots(a_r)_k}
{(b_1)_k\cdots(b_s)_k}\cdot \f{z^k}{k!}.$$
 After
reading the author's preprint (arXiv:1101.1050) involving
(1.4)-(1.6), in the email to the author on January 11, 2011 Wadim
Zudilin wrote: ``It's probably worth mentioning that the proofs of
Theorems 2.1, 2.4 and possibly Theorem 3.1 from your arXiv:1101.1050
can be simplified. Your congruences assume the form
 $$\left( \sum_{k=0}^{p-1} \frac{(a)_k (1-a)_k}{k!^2} x^k \right)^2
\e \sum_{k=0}^{p-1} \frac{(a)_k (1-a)_k}{k!^2} \binom{2k}{k}
(x(1-x))^k  \mod {p^2}, $$
 where $a=1/3$, $1/4$ or $1/6$ ($1/2$ is possible as well). Note that
  this follows from the identity
 $$ \left( \sum_{k=0}^\infty \frac{(a)_k (1-a)_k}{k!^2} x^k \right)^2
= \sum_{k=0}^\infty \frac{(a)_k (1-a)_k}{k!^2} \binom{2k}{k}
(x(1-x))^k $$
 truncated to $p$ terms. If $k>p/2$, then the coefficients on
 both sides
 are $0\mod {p^2}$ (as you already used in the proof, e.g., of Theorem
  2.1)." In the email on January 15, 2011 Zudilin wrote: The
  identity is a combination of the Gauss quadratic
transformation ([Ba, p.88, Eq. (2)])
$${}_2F_1\biggl(\begin{matrix} A, \, B \\ A+B+\frac12 \end{matrix}
\biggm|4z(1-z)\biggr) ={}_2F_1\biggl(\begin{matrix} 2A, \, 2B \\
A+B+\frac12 \end{matrix} \biggm|z\biggr),$$ and Clausen's original
identity ([Sl, p.75, Eq.(2.5.7)])
$${}_2F_1\biggl(\begin{matrix} a, \, b \\ a+b+\frac12 \end{matrix}
\biggm|z\biggr)^2 ={}_3F_2\biggl(\begin{matrix} 2a, \, 2b, \, a+b \\
a+b+\frac12, \, 2a+2b \end{matrix} \biggm|z\biggr).$$
 \par Let $p$ be an odd prime and $a\in\Bbb Z_p$. In this note,
  inspired by Zudilin's
comments we give an elementary proof of the following generalization
of (1.4)-(1.6):
$$\Big(\sum_{k=0}^{p-1}\b ak\b {-1-a}kx^k\Big)^2\e
\sum_{k=0}^{p-1}\b{2k}k\b ak\b {-1-a}k(x(1-x))^k\mod{p^2}.$$ As an
application, for $m\in\Bbb Z_p$ with
 $m\not\e 0\mod p$ we show that
 $$\sum_{k=0}^{p-1}\f{\b{2k}k\b ak\b{-1-a}k}{m^k}\e 0\mod p
 \qtq{implies}\sum_{k=0}^{p-1}\f{\b{2k}k\b ak\b{-1-a}k}{m^k}\e 0\mod
 {p^2}.$$

\section*{2. Main results}
\par\q\pro{Lemma 2.1 ([S3, Lemma 3.2])}  For any nonnegative integer $n$
we have
$$P_n(\sqrt{1+4x})^2=\sum_{k=0}^n\b nk\b {n+k}k\b{2k}kx^k.$$
\endpro
\par As Zudilin noted, Lemma 2.1 can also be deduced from
 Clausen's identity, Gauss' quadratic transformation for
hypergeometic series and (1.1).

 \pro{Theorem 2.1} Let $p$ be an odd prime and
$a\in\Bbb Z_p$. Then
$$\sum_{k=0}^{p-1}\b{2k}k\b ak\b{-1-a}kx^k
\e P_{\ap}(\sqrt{1-4x})^2\e P_{p-1-\ap}(\sqrt{1-4x})^2\mod p.$$
\endpro
Proof. By Lemma 2.1,
$$\align P_{\ap}(\sqrt{1-4x})^2
&=\sum_{k=0}^{\ap}\b{2k}k\b{\ap}k\b{\ap+k}k (-x)^k
\\&=\sum_{k=0}^{p-1}\b{2k}k\b{\ap}k\b{\ap+k}k (-x)^k
\\&\e \sum_{k=0}^{p-1}\b{2k}k\b ak\b{a+k}k (-x)^k
\\&=\sum_{k=0}^{p-1}\b{2k}k\b ak\b{-1-a}k x^k\mod p.
\endalign$$
To complete the proof, we note that $P_{p-1-\ap}(t)\e P_{\ap}(t)
\mod p$ by [S4, Lemma 2.2].
\par{\bf Remark 2.1} In the cases $a=-\f 12,-\f 13,-\f 14,-\f 16$,
Theorem 2.2 was given by the author in [S2, Theorem 4.1], [S3,
Theorems 3.2 and 4.2] and [S5, Theorem 4.1].

\pro{Corollary 2.1} Let $p$ be an odd prime and $a,x\in\Bbb Z_p$.
Then $\sum_{k=0}^{p-1}\b{2k}k\b ak\b{-1-a}kx^k\e 0\mod p$ implies
$$P_{\ap}(\sqrt{1-4x})\e P_{p-1-\ap}(\sqrt{1-4x}) \e 0\mod p.$$
\endpro
Proof. By (1.1),
$$P_n(\sqrt{1-4x})=\f 1{2^n}(\sqrt{1-4x})^{n-2[\f n2]}
\sum_{k=0}^{[n/2]}\b nk(-1)^k\b{2n-2k}n(1-4x)^{[\f n2]-k}.$$ Thus
$P_n(\sqrt{1-4x})^2\e 0\mod p$ implies $P_n(\sqrt{1-4x})\e 0\mod p$.
Now applying Theorem 2.1 we deduce the result.

\pro{Lemma 2.2} For any nonnegative integer $n$ we have
$$\sum_{k=0}^n\b ak\b{-1-a}k\b a{n-k}\b{-1-a}{n-k}
=\sum_{k=0}^n\b{2k}k\b ak\b{-1-a}k\b k{n-k}(-1)^{n-k}.$$
\endpro
Proof. Let $S_1(n)$ and $S_2(n)$ be the sums on the left and right
hands of the identity, respectively. Using Maple and the Zeilberger
algorithm we find that  for $i=1,2$,
$$n^3S_i(n)=(2n-1)(n^2-n-2a(a+1))S_i(n-1)+(n-1)(2a+n)(2a+2-n)S_i(n-2)\
(n\ge 2).$$ Since $S_1(0)=1=S_2(0)$ and $S_1(1)=-2a(a+1)=S_2(1)$,
applying the above we deduce that $S_1(n)=S_2(n)$.
\par{\bf Remark 2.2} In the cases $a=-\f 13,-\f 14,-\f 16$, the
identity was given by the author in [S2,S3,S5].

 \pro{Theorem 2.2} Let $p$ be an odd prime and $a\in\Bbb Z_p$. Then
$$\Big(\sum_{k=0}^{p-1}\b ak\b {-1-a}kx^k\Big)^2\e
\sum_{k=0}^{p-1}\b{2k}k\b ak\b {-1-a}k(x(1-x))^k\mod{p^2}.$$
\endpro
\par Taking $a=-\f 13,-\f 14,-\f 16$ in Theorem 2.2 and
then applying (1.7) we get (1.4)-(1.6). Thus, Theorem 2.2 is a
generalization of (1.4)-(1.6). \par Proof of Theorem 2.2. For
$k\in\{\f{p+1}2,\ldots,p-1\}$ we see that $p\mid \b{2k}k$ and
$$\align \b ak\b {-1-a}k&=(-1)^k\b ak\b{a+k}k
\\&=(-1)^k\f{(a+k)(a+k-1)\cdots (a-k+1)}{k!}\e 0\mod p.\endalign$$
Thus, applying Lemma 2.2 we see that
$$\align&\sum_{k=0}^{p-1}\b{2k}k\b ak\b {-1-a}k(x(1-x))^k
\\&\e \sum_{k=0}^{(p-1)/2}\b{2k}k\b ak\b {-1-a}k(x(1-x))^k
\\&=\sum_{k=0}^{(p-1)/2}\b{2k}k\b ak
\b {-1-a}kx^k\sum_{r=0}^k\b kr(-x)^r
\\&=\sum_{n=0}^{p-1}x^n
\sum_{k=0}^{\min\{n,\f{p-1}2\}}\b{2k}k\b ak\b{-1-a}k\b
k{n-k}(-1)^{n-k}
\\&\e\sum_{n=0}^{p-1}x^n
\sum_{k=0}^n\b{2k}k\b ak\b{-1-a}k\b k{n-k}(-1)^{n-k}
\\&=\sum_{n=0}^{p-1}x^n
\sum_{k=0}^n\b ak\b{-1-a}k\b a{n-k}\b{-1-a}{n-k}
\\&=\sum_{k=0}^{p-1}\b ak\b{-1-a}kx^k\sum_{n=k}^{p-1}
\b a{n-k}\b{-1-a}{n-k}x^{n-k}
\\&=\Big(\sum_{k=0}^{p-1}\b ak\b{-1-a}kx^k\Big)
\Big(\sum_{r=0}^{p-1} \b ar\b{-1-a}rx^r-\sum_{r=p-k}^{p-1} \b
ar\b{-1-a}rx^r\Big)\mod {p^2}.
\endalign$$
For $0\le k\le p-1$ we see that $p\mid \b ak$ for $k>\ap$. For $0\le
k\le \ap$ and $p-k\le r\le p-1$ we see that $r\ge p-k>p-1-\ap$ and
so $p\mid \b{-1-a}r$. Hence for $0\le k\le p-1$ and $p-k\le r\le
p-1$ we have $p\mid \b ak\b{-1-a}r$, $p\mid \b{-1-a}k\b ar$ and so
$p^2\mid \b ak\b{-1-a}k\b ar\b{-1-a}r$. Now combining all the above
we deduce the result.

 \pro{Theorem 2.3} Let $p$ be an odd prime and $m\in\Bbb Z_p$ with
 $m\not\e 0\mod p$. Then
 $$\sum_{k=0}^{p-1}\f{\b{2k}k\b ak\b{-1-a}k}{m^k}\e 0\mod p
 \qtq{implies}\sum_{k=0}^{p-1}\f{\b{2k}k\b ak\b{-1-a}k}{m^k}\e 0\mod
 {p^2}.$$
 \endpro
Proof. Suppose $\sum_{k=0}^{p-1}\f{\b{2k}k\b ak\b{-1-a}k}{m^k}\e
0\mod p$. By Corollary 2.1, $P_{\ap}(\sqrt{1-\f 4m})\e 0\mod p$.
Thus,
$$\align &\sum_{k=0}^{p-1}\b ak\b{-1-a}k\Ls{1-\sqrt{4/m}}2^k
\\&\e \sum_{k=0}^{p-1}\b {\ap}k\b{-1-\ap}k\Ls{1-\sqrt{4/m}}2^k
\\&=\sum_{k=0}^{\ap}\b{\ap}k\b{\ap+k}k\Ls{\sqrt{4/m}-1}2^k
\\&=P_{\ap}\Big(\sqrt{1-\f 4m}\Big)\e 0\mod p.
\endalign$$
By Theorem 2.2,
$$\Big(\sum_{k=0}^{p-1}\b ak\b{-1-a}k\Ls{1-\sqrt{4/m}}2^k\Big)^2\e
\sum_{k=0}^{p-1}\f{\b{2k}k\b ak\b{-1-a}k}{m^k}\mod {p^2}.$$ Thus the
result follows.

Taking $a=-\f 12,-\f 13,-\f 14,-\f 16$ in Theorem 2.3 and then
applying (1.7) we deduce  the following result.
 \pro{Corollary 2.2}
Let $p$ be an odd prime, and $m\in\Bbb Z_p$ with $m\not\e 0\mod p$.
Then
$$\align &\sum_{k=0}^{p-1}\f{\b{2k}k^3}{m^k}\e 0\mod p
\qtq{implies} \sum_{k=0}^{p-1}\f{\b{2k}k^3}{m^k}\e 0\mod {p^2},
\\&\sum_{k=0}^{p-1}\f{\b{2k}k^2\b{3k}k}{m^k}\e 0\mod p
\qtq{implies} \sum_{k=0}^{p-1} \f{\b{2k}k^2\b{3k}k}{m^k}\e 0\mod
{p^2},
\\&\sum_{k=0}^{p-1}\f{\b{2k}k^2\b{4k}{2k}}{m^k}\e 0\mod p
\qtq{implies} \sum_{k=0}^{p-1} \f{\b{2k}k^2\b{4k}{2k}}{m^k}\e 0\mod
{p^2}\endalign$$ and
$$\sum_{k=0}^{p-1}\f{\b{2k}k\b{3k}k\b{6k}{3k}}{m^k}\e 0\mod p
\qtq{implies} \sum_{k=0}^{p-1}\f{\b{2k}k\b{3k}k\b{6k}{3k}}{m^k}\e
0\mod {p^2}.$$
\endpro
\par We remark that Corollary 2.2 can be easily deduced from [S2, Theorem
4.1], [S3, Theorems 3.2 and 4.2] and [S5, Theorem 4.1].
 \pro{Theorem
2.4} Let $p$ be an odd prime and $u\in\Bbb Z_p$.
\par $(\t{\rm i})$ If $u\not\e \f 14,\f 1{16}\mod p$ and
$\sum_{k=0}^{p-1}\b{2k}k^2\b{3k}k\sls{u^2}{(1-4u)^3}^k\e 0\mod p$,
then
 $$\sum_{k=0}^{p-1}\b{2k}k^2\b{3k}k\Ls{-u}{(1-16u)^3}^k \e
0\mod{p^2}.$$
\par $(\t{\rm ii})$ If
$u\not\e -\f 13,-\f 1{27}\mod p$ and
$\sum_{k=0}^{p-1}\b{2k}k^2\b{4k}{2k}\sls {u^3}{(1+3u)^4}^k\e 0\mod
p$, then
$$\sum_{k=0}^{p-1}\b{2k}k^2\b{4k}{2k}\Ls u{(1+27u)^4}^k\e
0\mod{p^2}.$$
\endpro
Proof. This is immediate from Corollary 2.2 and [S7, Corollaries 3.1
and 5.2].

\pro{Corollary 2.3} Let $p>3$ be a prime. Then
$$\sum_{k=0}^{p-1}\f{\b{2k}k^2\b{3k}k}{1458^k}\e 0\mod {p^2}
\qtq{for}p\e 5\mod 6$$ and
$$\sum_{k=0}^{p-1}\f{\b{2k}k^2\b{3k}k}{15^{3k}}\e 0\mod {p^2}
\qtq{for}p\e 11,14\mod {15}.$$
\endpro
Proof. Taking $u=-\f 12$ in Theorem 2.4(i) and then applying (1.2)
we obtain the first congruence. Taking $u=1$ in Theorem 2.4(i) and
then applying [S2, Theorems 4.6] we obtain the second congruence.

{\bf Remark 2.3} Let $p$ be a prime of the form $6k+5$. In [S7] the
author proved that $\sum_{k=0}^{p-1}\f{\b{2k}k^2\b{3k}k}{1458^k}\e
0\mod p$ and conjectured that
$\sum_{k=0}^{p-1}\f{\b{2k}k^2\b{3k}k}{1458^k}\e 0\mod {p^3}$.
\par\q
\newline{\bf Acknowledgment}
\par\q
\newline The author is indebted to Wadim Zudilin for his
contribution to Theorem 2.2. See comments in Section 1.

\end{document}